\documentclass{article}
\pdfoutput=1
\usepackage{amssymb}
\usepackage{appendix}
\usepackage{latexsym}
\usepackage{mathrsfs}
\usepackage{enumerate}
\usepackage{extarrows}
\usepackage{graphicx}
\usepackage{caption}
\usepackage{supertabular}
\usepackage{subfigure}
\usepackage{setspace}
\usepackage{epstopdf}

\let\cal=\mathcal
\newtheorem{theo}{Theorem}[section]

\newtheorem{remark}[theo]{Remark}

\newtheorem{lemma}[theo]{Lemma}
\newtheorem{claim}[theo]{Claim}
\newtheorem{coro}[theo]{Corollary}

\newtheorem{fact}[theo]{Fact}
\newtheorem{defi}[theo]{Definition}

\def\q{\hspace*{\fill}$\Box$\medskip}
\def\defeq{\stackrel{\rm\scriptscriptstyle def}{=}}

\def\endproofbox{\hskip 1.3em\hfill\rule{6pt}{6pt}}

{%
\noindent{\it Proof}
}%
{%
 \quad\hfill\endproofbox\vspace*{2ex}
}
\begin{document}

\title{Non-jumping Tur\'an densities of hypergraphs }
\author{Zilong Yan \thanks{School of Mathematics, Hunan University, Changsha 410082, P.R. China. Email: zilongyan@hnu.edu.cn.} \and Yuejian Peng \thanks{ Corresponding author. School of Mathematics, Hunan University, Changsha, 410082, P.R. China. Email: ypeng1@hnu.edu.cn. \ Supported in part by National Natural Science Foundation of China (No. 11931002).}
}

\maketitle

\begin{abstract}
A real number $\alpha\in [0, 1)$ is a jump for an integer $r\ge 2$ if  there exists $c>0$ such that no number in  $(\alpha ,  \alpha + c)$ can be the Tur\'an density of a family of $r$-uniform graphs. A classical result of Erd\H os and Stone \cite{ES} implies that that every number in $[0, 1)$ is a jump for $r=2$. Erd\H os \cite{E64} also showed that every number in $[0, r!/r^r)$ is a jump for $r\ge 3$ and asked whether every number in $[0, 1)$ is a jump for $r\ge 3$.  Frankl and R\"odl \cite{FR84} gave a negative answer by showing a sequence of non-jumps  for every $r\ge 3$. After this, Erd\H os modified the question to be whether $\frac{r!}{r^r}$ is a jump for $r\ge 3$? What's the smallest non-jump? Frankl, Peng,  R\"odl and Talbot  \cite{FPRT} showed that ${5r!\over 2r^r}$ is a non-jump for $r\ge 3$.  Baber and Talbot \cite{BT0} showed that every $\alpha\in[0.2299, 0.2316)\cup [0.2871, \frac{8}{27})$ is a jump for $r=3$. Pikhurko \cite{Pikhurko2} showed that the set of all possible Tur\'an densities of $r$-uniform graphs has cardinality of the continuum for $r\ge 3$. However, whether $\frac{r!}{r^r}$ is a jump for $r\ge 3$ remains open, and $\frac{5r!}{2r^r}$ has remained the known smallest non-jump for $r\ge 3$. In this paper, we give a smaller non-jump by showing that ${54r!\over 25r^r}$ is a non-jump for $r\ge 3$. Furthermore, we give infinitely many irrational non-jumps for every $r\ge 3$.

\medskip

{\em Keywords}: Jumping number, Tur\'an density, hypergraph
\end{abstract}

\section{Introduction}
For a finite set $V$ and a positive integer $r$ we denote by ${{V \choose r}}$ the family of all $r$-subsets of $V$. An {\em $r$-uniform graph} ({\em $r$-graph}) $G$ is a set $V(G)$ of vertices together with a set $E(G) \subseteq {V(G) \choose r}$ of edges. The {\em density} of $G$ is defined by $d(G) = {\vert E(G)\vert\over \vert{V(G) \choose r}\vert}$. 
For a family ${\cal F}$ of $r$-graphs, an $r$-graph $G$ is called $\cal F$-free if it does not contain an isomorphic copy of any $r$-graph of $\cal F$. For a fixed positive integer $n$ and a family of $r$-graphs $\cal F$, the {\em Tur\'an number} of $\cal F$, denoted by $ex(n,\cal F)$, is the maximum number of edges in an $\cal F$-free $r$-graph on $n$ vertices. An averaging argument in \cite{KNS} by Katona, Nemetz, and Simonovits shows that the sequence ${ex(n, \cal F)\over {n\choose r}}$ is non-increasing. Hence $\lim_{n\to\infty}{ex(n, \cal F)\over {n\choose r}}$ exists. The {\em Tur\'{a}n density} of $\cal F$ is defined as $$\pi(\cal F)=\lim_{n\rightarrow\infty}{ex(n, \cal F) \over {n \choose r }}.$$
If $\cal F$ consists of a single $r$-graph $F$, we simply write $ex(n, \{F\})$ and $\pi(\{F\})$ as $ex(n, F)$ and $\pi(F)$. Denote
$$\Pi_{\infty}^{r}=\{ \pi(\cal F): \cal F {\rm \ is \ a \ family \ of \ } r{\rm-uniform \ graphs} \}, $$
$$\Pi_{fin}^{r}=\{\pi(\cal F): \cal F {\rm \ is \ a \ finite \ family \ of \ } r{\rm-uniform \ graphs}, \}$$
and
$$\Pi_{t}^{r}=\{ \pi(\cal F): \cal F {\rm \ is \ a \ family \ of \ } r{\rm-uniform \ graphs \ and }\ \vert \cal F \vert\le t  \}. $$
Clearly,
$$\Pi_{1}^{r}\subseteq \Pi_{2}^{r}\subseteq \cdots \subseteq\Pi_{fin}^{r}\subseteq \Pi_{\infty}^{r}.$$
Finding good estimation for Tur\'an densities in hypergraphs ($r\ge 3$) is believed to be one of the most challenging problems in extremal combinatorics. The following concept concerns the accumulation points of the set $\Pi_{\infty}^{r}$.
\begin{defi}
A real number $\alpha\in [0, 1)$ is a {\em jump} for an integer $r\ge 2$ if there exists a constant $c>0$ such that for any $\epsilon>0$ and any integer $m \ge r$, there exists an integer $n_0(\epsilon, m)$ such that any $r$-uniform graph with $n > n_0(\epsilon, m)$ vertices and density $\ge \alpha + \epsilon$ contains a subgraph with $m$ vertices and density $\ge \alpha + c$.
\end{defi}

This concept describes where the set of `jumps' is closely related to Tur\'an densities. It was shown in \cite{FR84} that  $\alpha$ is a jump for $r$ if and only if there exists $c>0$ such that $\Pi_{\infty}^r \cap (\alpha, \alpha+c)=\emptyset$. So every non-jump is an accumulation point of $\Pi_{\infty}^{r}$.
For 2-graphs, Erd\H{o}s-Stone-Simonovits \cite{ESi,ES} determined the Tur\'an numbers of all non-bipartite graphs asymptotically. Their result implies that $$\Pi_{\infty}^{2}=\Pi_{fin}^{2}=\Pi_{1}^{2}=\{0, {1 \over 2}, {2 \over 3}, ..., {l-1 \over l}, ...\}.$$
This implies that every $\alpha\in [ 0,  1)$ is a jump for $r=2$. For $r\geq 3$, Erd\H{o}s \cite{E64} proved that every $\alpha\in [0, r!/r^r)$ is a jump. Furthermore, Erd\H{o}s proposed the {\it jumping constant conjecture}: Every $\alpha\in [0, 1)$ is a jump for every integer $r \geq 2$. In \cite{FR84}, Frankl and R\"{o}dl disproved the Conjecture by showing that $\displaystyle{1-\frac{1}{l^{r-1}}}$ is not a jump for $r$ if $r\ge 3$ and $l>2r$. However, there are still a lot of unknowns on whether a number is a jump for $r\ge 3$. A well-known open question of Erd\H{o}s is whether $r!/r^r$ is a jump for $r\ge 3$ and what is the smallest non-jump? Another question raised in \cite{FPRT} is whether there is an interval  of non-jumps for some $r\ge 3$? Both questions seem to be very challenging. Frankl-Peng-R\"{o}dl-Talbot \cite{FPRT} showed that ${5r! \over 2 r^r}$ is a non-jump for $r\ge 3$. Baber and Talbot \cite{BT0} showed that for $r=3$ every $\alpha\in[0.2299, 0.2316)\cup [0.2871, \frac{8}{27})$ is a jump. Pikhurko \cite{Pikhurko2} showed that $\Pi_{\infty}^r$ has cardinality of the continuum for $r\ge 3$. However, whether $\frac{r!}{r^r}$ is a jump remains open. Regarding the first question, we determine a non-jump smaller than ${5r! \over 2 r^r}$ for $r\ge 3$.

\begin{theo}\label{theo}
$\frac{12}{25}$ is not a jump for $r=3$.
\end{theo}

In \cite{jumpgeneral}, a way to generate non-jumps for every $p\ge r$ based on a non-jump for $r$ was given. The following result was shown there.

\begin{theo}\label{resultg}\cite{jumpgeneral}
Let $p\ge r\ge 3$ be positive integers. If $\alpha\cdot {r! \over r^r} $ is a non-jump  for $r$, then $\alpha \cdot {p! \over p^p}$ is a non-jump for $p$.
\end{theo}

Combining Theorems \ref{theo} and \ref{resultg}, we will get

\begin{coro}
${54r! \over 25r^r}$ is a non-jump for $r\ge 3$.
\end{coro}


Chung and Graham \cite{FG} proposed the conjecture that every element in $\Pi_{fin}^{r}$ is a rational number. Baber and Talbot \cite{BT}, and Pikhurko \cite{Pikhurko2} disproved this conjecture independently by showing that there is an irrational number in $\Pi_{fin}^{r}$. Baber and Talbot asked that whether there is an irrational number in $\Pi_1^r$. Recently, Yan and Peng \cite{YP} showed that there is an irrational number in $\Pi_1^3$ and Wu-Peng \cite{WP} showed that there is an irrational number in $\Pi_1^4$. Pikhurko \cite{Pikhurko2} showed that $\Pi_{\infty}^r$ is closed which implies that every non-jump is a Tur\'an density (a Tur\'an density may not be a non-jump). Brown and Simonovits \cite{BS} showed that the Lagrangian of an $r$-uniform hypergraph is in $\Pi_{\infty}^r$ also indicating the existence of irrational numbers in $\Pi_{\infty}^r$. No irrational non-jump has been previously given. In this paper, we will give a infinite sequence of irrational non-jumps for $r=3$.

\begin{theo}\label{theo1}
Let $k\ge 2$ be an integer. Then $\alpha_k=\frac{2k-6k^3+4k^4-k\sqrt{4k - 1}+4k^2\sqrt{4k - 1}}{(2k^2+1)^2}$ is not a jump for $r=3$.
\end{theo}

Combining Theorem \ref{theo1} and Theorem \ref{resultg}, we can also get corresponding non-jumps for $r\ge 3$.

The proof of Theorem \ref{theo} and \ref{theo1} will be given in Section \ref{prooftheo} and Section \ref{prooftheo1}, respectively. Both proofs applied an approach developed by Frankl and R\"odl in \cite{FR84}. The crucial part in our proof is to give a `proper' construction. In the following section, we will introduce some preliminary results and sketch the idea of the proof.

\section{Preliminaries and Sketch of the proof}
\subsection{Karush-Kuhn-Tucker Conditions}
Let us consider the optimisation problem:
\begin{flushleft}
\quad\quad\quad\quad\quad maximise\ $f(x)$\\
\quad\quad\quad\quad\quad subject to $g_i(x)\leq 0$, $i=1,\dots,m,$\quad\quad\quad\quad\quad\quad\quad\quad\quad\quad\quad\quad\quad (3.1)
\end{flushleft}
where $x\in \mathbb{R}^n$ and $f$ and $g_i$ are differentiable functions from $\mathbb{R}^n$ to $\mathbb{R}$ for all $i$. Let $\nabla{f(x)}$ be the gradient of $f$ at $x$ i.e., the vector in $\mathbb{R}^n$ whose $i$th coordinate is ${\partial\over\partial{x_i}}f(x)$. We say that KKT conditions hold at $x^*\in \mathbb{R}^n$ if there exist $\lambda_1,\dots\lambda_m\in \mathbb{R}$ such that
\begin{enumerate}
\item $\nabla{f(x^*)}=\sum_{i=1}^m\lambda_i\nabla{g_i(x^*)},$
\item $\lambda_i\geq 0$ for $i=1,\dots,m,$
\item $\lambda_ig_i(x^*)=0$ for $i=1,\dots,m.$
\end{enumerate}
We call the constraints linear if $g_1,\dots,g_m$ are all affine functions.
\begin{theo}\label{KKT}(\cite{BV},\cite{Jenssen})
If the constraints of (3.1) are linear, then any optimal solution to (3.1) must satisfy the KKT conditions.
\end{theo}

\subsection{Properties of the Lagrangian function}
In this section we will give the definition of the Lagrangian of an $r$-uniform graph, which is a helpful tool in our proof.

\begin{defi}
For an $r$-uniform graph $G$ with vertex set $\{1,2,\ldots,n\}$, edge set $E(G)$ and a vector $\vec{x}=(x_1,\ldots,x_n) \in R^n$, define the Lagrangian fuction
$$\lambda (G,\vec{x})= \sum_{\{i_1,\ldots,i_r\}\in E(G)}x_{i_1}x_{i_2}\ldots x_{i_r}.$$
\end{defi}
Let $S=\{\vec{x}=(x_1,x_2,\ldots ,x_n): \sum_{i=1}^{n} x_i =1, x_i \ge 0 {\rm \ for \ } i=1,2,\ldots , n  \}$. The Lagrangian of $G$, denoted by $\lambda (G)$, is defined as
$$\lambda (G) = \max \{\lambda (G, \vec{x}): \vec{x} \in S \}.$$
A vector $\vec{x}\in S$ is called a {\em feasible vector} on $G$, and $x_i$ is called the {\em weight} of the vertex $i$.
 A feasible vector is called {\em optimal} if $\lambda (G, \vec{y})=\lambda(G)$.

\begin{fact}\label{mono}
If $G_1\subseteq G_2$, then $$\lambda (G_1) \le \lambda (G_2).$$
\end{fact}

\begin{fact}(\cite{FR84})\label{fact2}
Let $G$ be an $r$-graph on $[n]$. Let $\vec{x}=(x_1,x_2,\dots,x_n)$ be an optimal vector on $G$. Then
$$ \frac{\partial \lambda (G, \vec{x})}{\partial x_i}=r\lambda(G)$$
for every $i \in [n]$ satisfying $x_i>0$.
\end{fact}

Given an $r$-graph $G$, and $i, j\in V(G),$ define $$L_G(j\setminus i)=\{e: i\notin e, e\cup\{j\}\in E(G)\:and\: e\cup\{i\}\notin E(G)\}.$$

\begin{fact}\label{symmetry}(\cite{FR84})
Let $G$ be an $r$-graph on $[n]$. Let $\vec{x}=(x_1,x_2,\dots,x_n)$ be a feasible vector on $G$, and $i,j\in [n]$, $i\neq j$ satisfy $L_G(i \setminus j)=L_G(j \setminus i)=\emptyset$. Let $\vec{y}=(y_1,y_2,\dots,y_n)$ be defined by letting $y_\ell=x_\ell$ for every $\ell \in [n]\setminus \{i,j\}$ and $y_i=y_j={1 \over 2}(x_i+x_j)$.
Then $\lambda(G,\vec{y})\geq \lambda(G,\vec{x})$. Furthermore, if the pair $\{i,j\}$ is contained in an edge of $G$, $x_i>0$ for each $1\le i\le n$,  and $\lambda(G,\vec{y})=\lambda(G,\vec{x})$, then $x_i=x_j$.
\end{fact}

We also note that for an $r$-graph $G$ with $n$ vertices, if we take $\vec{u}=(u_1, \ldots, u_n)$, where each $u_i={1\over n}$, then
$$\lambda(G)\ge \lambda(G, \vec{u})={\vert E(G)\vert \over n^r} \ge {d(G) \over r!}-\epsilon$$ for $n\ge n'(\epsilon)$, where $n'(\epsilon)$ is a sufficiently large integer. On the other hand, the blow-up of an $r$-uniform graph $G$ will allow us to construct  $r$-uniform graphs with  large  number of vertices and  density close to $r!\lambda (G)$.

\begin{defi}
Let $G$ be an $r$-uniform graph with $V(G) =\{1,2,\ldots,t\}$ and $(n_1, \ldots, n_t)$ be a positive integer vector. Define the $(n_1, \ldots, n_t)$ blow-up of $G$, $(n_1, \ldots, n_t)\otimes G$ as a $t$-partite $r$-uniform graph with vertex set $V_1\cup \ldots \cup V_t, |V_i|=n_i, 1\leq i\leq t$, and edge set $E((n_1, \ldots, n_t)\otimes G)=\{\{v_{i_1}, v_{i_2},\ldots, v_{i_r}\}, { \ \rm where \ } \{i_1,i_2,\ldots, i_r\} \in E(G) {\ \rm and \ } v_{i_k} \in V_{i_k} {\rm \ for} \ 1\le k\le r \}$.
\end{defi}

\begin{remark} (\cite{FR84})\label{remarkblow}
Let $G$ be an $r$-uniform graph with $t$ vertices and $\vec{y}=(y_1, \ldots, y_t)$ be an optimal vector for $\lambda(G)$. Then for any $\epsilon >0$, there exists an integer $n_1(\epsilon)$, such that for any integer $n\ge n_1(\epsilon)$,
\begin{equation}\label{blowden}
d((\lceil ny_1\rceil, \lceil ny_2\rceil, \ldots, \lceil ny_t\rceil)\otimes G)\ge r!\lambda(G)-\epsilon.
\end{equation}
\end{remark}

Let us also state a fact which follows directly from the definition of the Lagrangian.

\begin{fact}(\cite{FR84})\label{lblow}
For every $r$-uniform graph $G$ and every positive integer $n$, $\lambda((n, n, \ldots,n)\otimes G) =\lambda (G)$ holds.
\end{fact}

Lemma \ref{arrow} in \cite{FR84} gives a necessary and sufficient condition for a number $\alpha$ to be a jump.

\begin{lemma}\label{arrow}(\cite{FR84})
The following two properties are equivalent.
\begin{enumerate}
\item $\alpha$ is a jump for $r$.
\item There exists  some finite family $\cal F$ of $r$-uniform graphs satisfying $\pi({\cal F})\le \alpha$ and $\displaystyle{\lambda (F)> \frac{\alpha}{r!}}$ for all $F \in \cal F$.
\end{enumerate}
\end{lemma}

\subsection{Sketch of the proofs of Theorem \ref{theo} and \ref{theo1}}
The general approach in proving Theorem \ref{theo} and Theorem  \ref{theo1} is sketched as follows: Let $\alpha$ be a number to be proved to be a non-jump for $r=3$. Assuming that $\alpha$ is a jump for $r=3$, we will derive a contradiction by the following steps.

Step1. Construct a `proper' $3$-uniform graph
$G^*(t)$ with the Lagrangian at least ${\alpha \over 6}+\epsilon$ for some $\epsilon>0$. Then we `blow up' it to a $3$-uniform graph, say $\vec{m}\otimes G^*(t)$ with large enough number of vertices and density $\ge \alpha+\epsilon$ (see Remark \ref{remarkblow}). If $\alpha$ is a jump, then by Lemma \ref{arrow}, there exists some finite family ${\mathcal F}$ of $3$-uniform graphs with Lagrangians $>{\alpha \over 6}$ and $\pi(\cal{F})\le \alpha$. So $\vec{m}\otimes G^*(t)$ must contain some member of ${\mathcal F}$ as a subgraph.

Step 2. We show that any subgraph of $G^*(t)$ with the number of vertices not greater than $\max\{\vert V(F)\vert, F \in {\mathcal F}\}$ has the Lagrangian $\le {\alpha \over 6}$ and derive a contradiction.

\bigskip

The crucial part is to construct an $r$-uniform graph satisfying the properties in both Steps 1 and 2. Generally, whenever we find such a construction, we can obtain a corresponding non-jump. This method was first developed by Frankl and R\"odl in \cite{FR84}. The technical part in the proof is to show that the construction satisfies the property in Step 2.

\section{Proof of Theorem \ref{theo}} \label{prooftheo}
{\em Proof.} Suppose that ${12\over 25}$ is a jump for $r=3$. By Lemma \ref{arrow}, there exists a finite collection $\cal F$ of $3$-uniform graphs satisfying the following:
\begin{enumerate}
\item[i)] $\displaystyle \lambda (F) > {2 \over 25} $ for all $F \in \cal F$, and
\item[ii)]  $\pi ({\cal F})\le {12\over 25}$.
\end{enumerate}

Let $G(t)=(V, E)$ be the 3-uniform defined as follows. The vertex set $V=V_1\cup V_2\cup V_{3}$, where $|V_1|=|V_2|=\frac{2t}{5}$ and $|V_{3}|=\frac{t}{5}$ and the value of $t$ will be determined later. The edge set of $G(t)$ is
$$\bigg(V_1 \times V_2 \times V_3\bigg) \bigcup \bigg( {V_1\choose 2}\times V_2\bigg) \bigcup \bigg({V_2\choose 2}\times V_3\bigg),$$
i.e., the edges consisting of one vertex from each $V_1, V_2$ and $V_3$, or two vertices from $V_1$ and one vertex from $V_2$, or two vertices from $V_2$ and one vertex from $V_3$. Then
\begin{eqnarray}\label{eg1}
\vert E(G(t))\vert&=&\frac{2t^3}{25}-\frac{3t^2}{25}.
\end{eqnarray}
We will apply the following lemma from \cite{FR84}.

\begin{lemma}\label{add}\cite{FR84}
For any $c\ge 0$ and any integer $s\ge r$, there exists $t_0(s, c)$ such that for every $t\ge t_0(s, c)$, there exists an $r$-uniform graph $A=A(s, c, t)$ satisfying:
\begin{enumerate}
\item $|V(A)|=t,$
\item $|E(A)|\geq ct^{r-1},$
\item For all $V_0 \subset V(A), r\leq |V_0| \leq s$, we have $|E(A)\cap {V_0 \choose r}| \leq |V_0|-r+1$.
\end{enumerate}
\end{lemma}

Set $s= {\rm max}_{F \in {\cal F}} |V(F)|$ and $c=1$. Let $r=3$ in Lemma \ref{add}, $t_0(s, 1)$ be given as in Lemma \ref{add} and $\frac{2t}{5}\ge t_0(s, 1)$. The $3$-uniform graph $G^*(t)$ is obtained by adding $A(s, 1, \frac{2t}{5})$ to the $3$-uniform hypergraph $G(t)$ in $V_1$. Then
$$\lambda(G^*(t))\ge\lambda(G^*(t), (\frac{1}{t}, \frac{1}{t}, \dots, \frac{1}{t}))=\frac{\big\vert E(G^*(t))\big\vert}{t^3}.$$
In view of the construction of $G^*(t)$ and equation (\ref{eg1}), we have
\begin{eqnarray}\label{egAl}
\frac{\big\vert E(G^*(t))\big\vert}{t^3}&\ge&\frac{\big\vert E(G(t))\big\vert}{t^3}+\bigg(\frac{2t}{5}\bigg)^2/{t^3} \nonumber \\
&\ge& \frac{2}{25}+\frac{1}{25t}.
\end{eqnarray}

Now suppose $\vec{y}=(y_1, y_2, ..., y_t)$ is an optimal vector of $\lambda(G^*(t))$. Let $n$ be large enough. By Remark \ref{remarkblow}, $3$-uniform graph $S_n=(\lfloor ny_1\rfloor, \ldots, \lfloor ny_{n}\rfloor)\otimes G^*(t)$ has density at least${12\over 25}+\frac{1}{50t}.$ Since $\pi({\cal F})\le {12\over 25}$, some member $F$ of $\cal F$ is a subgraph of $S_n$ for $n$ sufficiently large. For such $F\in \cal F$, there exists a subgraph $M$ of $G^*(t)$ with $|V(M)|\le |V(F)|\leq s$ so that $F\subset (s, s, \ldots, s) \otimes M$. By Fact \ref{mono} and Fact \ref{lblow}, we have
\begin{equation}\label{lambdasmall0}
   \lambda(F)\overset{Fact \ref{mono}}{\le}\lambda ((s, s, \ldots, s) \otimes M)\overset{Fact \ref{lblow}}{=} \lambda (M).
\end{equation}
The following lemma will be proved in Section \ref{prooflemmaresult01}.

\begin{lemma}\label{lemmaresult01}
 Let $M$ be any subgraph of $G^*(t)$ with $|V(M)| \leq s$. Then
\begin{equation}
\lambda (M) \leq \frac{2}{25}
\end{equation}
holds.
\end{lemma}

Assuming that Lemma \ref{lemmaresult01} is true and applying Lemma \ref{lemmaresult01} to (\ref{lambdasmall0}), we have $$\lambda(F) \le {2 \over 25}$$ which contradicts our choice of $F$, i.e., contradicts that $\displaystyle \lambda(F) >{2 \over 25}$ for all $F \in \cal F$.  \q

\bigskip

To complete the proof of Theorem \ref{theo1}, what remains is to show Lemma \ref{lemmaresult01}.

\subsection{Proof of Lemma \ref{lemmaresult01}}\label{prooflemmaresult01}
By Fact \ref{mono}, we may assume that $M$ is an induced subgraph of $G^*(t)$. Let $$U_i=V(M)\cap V_i=\{v_1^i, v_2^i, \cdots, v_{s_i}^i\}.$$ So $s=s_1+s_2+s_3$. Let $\vec{z}=(z_1, z_2, ..., z_s)$ be an optimal vector for $M$. Without loss of generality, assume that $v_1^1$, $v_2^1$, $\cdots$, $v_{s+2}^1$ have the $s+2$  largest weights. Then replacing the $s$ edges in $M[U_1]$ by $v_1^1v_2^1v_3^1, v_1^1v_2^1v_4^1, \dots, v_1^1v_2^1v_{s+2}^1$ doesn't decrease the Lagrangian. So we have the following claim similar to Claim 4.4 in \cite{FR84}.
\begin{claim}\label{reduce0}
If $N$ is the $3$-uniform graph formed from $M$ by removing the edges contained in $U_{1}$ and inserting the edges $v^1_{1}v^1_2v^1_{j}$, where $3\leq j \leq s_1$, then $\lambda(M)\leq \lambda(N)$.
\end{claim}

By Claim \ref{reduce0}, the proof of Lemma \ref{lemmaresult01} will be completed if we show that $\lambda(N)\leq {2 \over 25}$. By Lemma \ref{symmetry}, we can obtain an optimal vector $\vec{z}$ of $\lambda(N)$ such that
\begin{equation} \label{weights}
w(v_1^1)=w(v_2^1)\defeq\frac{a}{2}, \ \ w(v_3^1)=w(v_4^1)=\cdots =w(v^1_{s_1}) \defeq\frac{b}{s_1-2},
\end{equation}
where $w(v)$ denotes the component of $\vec{z}$  corresponding to vertex $v$.

Let $c$, $d$ be the sum of the components of $\vec{z}$ corresponding to all vertices in $U_2$ and $U_3$, respectively. Note that
$$a+b+c+d=1.$$
Then
\begin{eqnarray*}
\lambda(N)\le\bigg(\frac{a^2}{4}+ab+\frac{b^2}{2}\bigg)c+(a+b)cd+\frac{c^2d}{2}+\frac{a^2}{4}b=\lambda(a, b, c, d).
\end{eqnarray*}
From now on, we assume that $(a, b, c, d)$ is an optimal vector for $\lambda(a, b, c, d)$.

If $c=0$, then $$\lambda(a, b, c, d)=\frac{a^2b}{4}\leq \frac{1}{8}\bigg(\frac{a+a+2b}{3}\bigg)^3\leq \frac{1}{27}.$$
So we may assume that $c>0$.

If $a=0$, then $$\lambda(a, b, c, d)=\frac{b^2c}{2}+bcd+\frac{c^2d}{2}\triangleq\lambda.$$ If $b=0$, then $\lambda=\frac{c^2d}{2}\leq \frac{2}{27}.$ Similarly we have $d>0$. So we may assume that $b, c, d>0$ in this case. By Theorem \ref{KKT}, we have
$$\frac{\partial\lambda}{\partial b}=\frac{\partial\lambda}{\partial c}=\frac{\partial\lambda}{\partial d},$$ so
$$bc+cd=\frac{b^2}{2}+bd+cd=\frac{c^2}{2}+bc.$$
Combining with $b+c+d=1$, we have $b=c=2d=0.4$, and $\lambda=\frac{2}{25}$. So we may assume that $a>0$.

If $b=0$, then $$\lambda(a, b, c, d)=\frac{a^2c}{4}+acd+\frac{c^2d}{2}<\frac{a^2c}{2}+acd+\frac{c^2d}{2}\leq \frac{2}{25}$$
as we have shown that $\frac{b^2c}{2}+bcd+\frac{c^2d}{2}\le \frac{2}{25}$. So we may assume that $b>0$.

We will prove that $d>0$ next. If $d=0$, then $$\lambda(a, b, c, d)=\bigg(\frac{a^2}{4}+ab+\frac{b^2}{2}\bigg)c+\frac{a^2}{4}b\triangleq\lambda.$$
By Theorem \ref{KKT}, we have
$$\frac{\partial\lambda}{\partial a}=\frac{\partial\lambda}{\partial b}.$$ So
$$\bigg(\frac{a}{2}+b\bigg)c+\frac{ab}{2}=(a+b)c+\frac{a^2}{4},$$
i.e., $a=2b-2c$. Since $a+b+c=1$, then $c=3b-1$ and $a=2-4b$. So
\begin{eqnarray*}
\lambda&\le&\bigg(\frac{a^2}{4}+ab+\frac{b^2}{2}\bigg)c+\frac{a^2}{4}b \\
&=&\frac{11b^3}{2}-\frac{21b^2}{2}+6b-1=f(b). \\
f'(b)&=&\frac{33b^2}{2}-21b+6.
\end{eqnarray*}
Since $a, c>0$, then $\frac{1}{3}\leq b\leq \frac{1}{2}$. Therefore $f(b)$ is increasing in $[\frac{1}{3}, \frac{7-\sqrt5}{11}]$ and decreasing in $[\frac{7-\sqrt5}{11}, \frac{1}{2}]$. Then $\lambda<f(\frac{7-\sqrt5}{11})<0.076.$

So we assume that $a, b, c, d>0$, then we have
\begin{eqnarray*}
\frac{\partial\lambda(a, b, c, d)}{\partial a}&=&\bigg(\frac{a}{2}+b\bigg)c+cd+\frac{ab}{2}, \\
\frac{\partial\lambda(a, b, c, d)}{\partial b}&=&(a+b)c+cd+\frac{a^2}{4}, \\
\frac{\partial\lambda(a, b, c, d)}{\partial c}&=&\frac{a^2}{4}+ab+\frac{b^2}{2}+ad+bd+cd,  \\
\frac{\partial\lambda(a, b, c, d)}{\partial d}&=&ac+bc+\frac{c^2}{2}, \\
d&=&1-a-b-c.
\end{eqnarray*}
By Theorem \ref{KKT}, we have $$\frac{\partial\lambda(a, b, c, d)}{\partial a}=\frac{\partial\lambda(a, b, c, d)}{\partial b},$$
and we get $a=2b-2c$. Therefore $d=1-a-b-c=1-3b+c$.
By $$\frac{\partial\lambda(a, b, c, d)}{\partial b}=\frac{\partial\lambda(a, b, c, d)}{\partial d},$$
we get $\frac{c^2}{2}=cd+\frac{a^2}{4}=c-3bc+c^2+b^2-2bc+c^2,$ so $$c=\frac{5b-1\pm \sqrt{19b^2-10b+1}}{3}.$$
By $$\frac{\partial\lambda(a, b, c, d)}{\partial b}=\frac{\partial\lambda(a, b, c, d)}{\partial c},$$ we get $$c=\frac{13b^2-6b}{8b-4}.$$
Therefore, $$\frac{13b^2-6b}{8b-4}=\frac{5b-1\pm \sqrt{19b^2-10b+1}}{3}.$$
By direct calculation, we have  $(-b^2+10b-4)^2=\bigg(\pm(8b-4)\sqrt{19b^2-10b+1}\bigg)^2.$ Simplifying, we get $$9b(5b-2)(9b-4)(3b-2)=0.$$
If $b=\frac{2}{5}$, then $c=\frac{13b^2-6b}{8b-4}=\frac{2}{5}$ and $a=2b-2c=0$, a contradiction. \\
If $b=\frac{4}{9}$, then $c=\frac{13b^2-6b}{8b-4}=\frac{2}{9}$, $a=2b-2c=\frac{4}{9}$ and $d=-\frac{1}{9}$, a contradiction. \\
If $b=\frac{2}{3}$, then $c<\frac{1}{3}$ and $a=2b-2c>\frac{2}{3}$ and $d<0$, a contradiction.
\q

\bigskip

\section{Proof of Theorem \ref{theo1}} \label{prooftheo1}
Let $B(2k, n)$ be the 3-graph with vertex set $[n]$ and edge set $E(B(2k, n))=\{e\in {[n]\choose 3}: e\cap [2k]\not=\emptyset\}$. Let $\alpha_k=\frac{2k-6k^3+4k^4-k\sqrt{4k - 1}+4k^2\sqrt{4k - 1}}{(2k^2+1)^2}$. We first show that $\alpha_k=6\lim_{n\to\infty}\lambda(B(2k, n))$.

Let $\vec{x}=\{x_1,x_2,\dots,x_n\}$ be an optimal vector of $\lambda(B(2k, n))$. Let $x_1+x_2+\cdots+x_{2k}=a$ and $b=1-a$. Then
\begin{eqnarray*}
\lim_{n\to\infty}\lambda(B(2k, n))&=&\bigg(\frac{a}{2k}\bigg)^3{2k\choose 3}+\bigg(\frac{a}{2k}\bigg)^2{2k\choose 2}(1-a)+a\frac{(1-a)^2}{2}=f(a)\\
f'(a)&=&(\frac{1}{4k^2}+\frac{1}{2})a^2-(\frac{1}{2k}+1)a+\frac{1}{2}.
\end{eqnarray*}
Note that $f(a)$ is increasing in $[0, \frac{2k^2+k-k\sqrt{4k-1}}{2k^2+1}]$ and decreasing in $[\frac{2k^2+k-k\sqrt{4k-1}}{2k^2+1}, 1]$. Therefore
\begin{eqnarray*}
f(\frac{2k^2+k-k\sqrt{4k-1}}{2k^2+1})&=&\frac{2k-6k^3+4k^4-k\sqrt{4k - 1}+4k^2\sqrt{4k - 1}}{6(2k^2+1)^2}\\
&=&\frac{\alpha_k}{6}.
\end{eqnarray*}
Since $k\ge 1$ and $4k-1$ is not a square number (a square number is 0 or 1 mod(4)), then $\alpha_k$ is an irrational number.

{\em Proof of Theorem \ref{theo1}.} Suppose that $\alpha_k$ is a jump. By Lemma \ref{arrow},  there exists a finite collection $\cal F$ of $3$-uniform graphs satisfying the following:
\begin{enumerate}
\item[i)] $\displaystyle \lambda (F) > {\alpha_k \over 6} $ for all $F \in \cal F$, and
\item[ii)]  $\pi ({\cal F})\le \alpha_k$.
\end{enumerate}

Let $G(t)=(V, E)$ be the 3-uniform defined as follows. The vertex set $V=V_1\cup V_2\cdots\cup V_{2k}\cup V_{2k+1}$, where $|V_1|=|V_2|=\cdots=|V_{2k}|=\frac{2k+1-\sqrt{4k-1}}{4k^2+2}t$ and $|V_{2k+1}|=\frac{k\sqrt{4k-1}+1-k}{2k^2+1}t$. The edge set of $G(t)$ is
$$\bigcup_{1\le i_1<i_2<i_3\le 2k} (V_{i_1}\times V_{i_2} \times V_{i_3})\bigcup_{1\le i_1<i_2\le 2k} (V_{i_1}\times V_{i_2} \times V_{2k+1})\bigcup_{1\le i_1\le 2k} (V_{i_1}\times{V_{2k+1}\choose 2}).$$
Then
\begin{eqnarray}\label{egl}
\vert E(G(t))\vert&=&\frac{2k-6k^3+4k^4-k\sqrt{4k - 1}+4k^2\sqrt{4k - 1}}{6(2k^2+1)^2}t^3 \nonumber \\
&+&\frac{-3k-6k^2+18k^3+3k\sqrt{4k - 1}-6k^2\sqrt{4k - 1}-6k^3\sqrt{4k - 1}}{6(2k^2+1)^2}t^2. \nonumber
\end{eqnarray}
Let $\vec{u}=(u_1, \ldots,u_{t})$, where $u_i=\frac{1}{t}$ for $1\le i\le t$,
then
\begin{eqnarray}\label{egl}
\lambda(G(t))&\ge&\lambda(G(t),\vec{u})={\vert E(G)\vert \over t^3} \nonumber \\
&=&\frac{2k-6k^3+4k^4-k\sqrt{4k - 1}+4k^2\sqrt{4k - 1}}{6(2k^2+1)^2} \nonumber \\
&+&\frac{-3k-6k^2+18k^3+3k\sqrt{4k - 1}-6k^2\sqrt{4k - 1}-6k^3\sqrt{4k - 1}}{6(2k^2+1)^2t}\\
&=&\frac{\alpha_k}{6}-\frac{c_0}{t}, \nonumber
\end{eqnarray}
where $c_0=\frac{3k+6k^2-18k^3-3k\sqrt{4k - 1}+6k^2\sqrt{4k - 1}+6k^3\sqrt{4k - 1}}{6(2k^2+1)^2}>0$.

Set $s= {\rm max}_{F \in {\cal F}} |V(F)|$ and $c=k$. Let $r=3$ in Lemma \ref{add} and $t_0(s, k)$ be given as in Lemma \ref{add}. Take an integer $t>\frac{2k^2+1}{k\sqrt{4k-1}+1-k}t_0$. The $3$-uniform graph $G^*(t)$ is obtained by adding $A(s, k)$ to the $3$-uniform hypergraph $G(t)$ in $V_{2k+1}$. Then
$$\lambda(G^*(t)) \ge  \frac{\big\vert E(G^*(t))\big\vert}{t^3}.$$
In view of the construction of $G^*(t)$ and equation (\ref{egl}), we have
\begin{eqnarray}\label{egAl}
\frac{\big\vert E(G^*(t))\big\vert}{t^3}&\ge&\frac{\big\vert E(G(t))\big\vert}{t^3}+\frac{k(\frac{k\sqrt{4k-1}+1-k}{2k^2+1}t)^2}{t^3} \nonumber \\
&{(\ref{egl}) \atop =}&\frac{2k-6k^3+4k^4-k\sqrt{4k - 1}+4k^2\sqrt{4k - 1}}{6(2k^2+1)^2} \nonumber \\
&+&\frac{3k-18k^2+18k^3+24k^4+3k\sqrt{4k - 1}+6k^2\sqrt{4k - 1}-18k^3\sqrt{4k - 1}}{6(2k^2+1)^2t} \nonumber \\
&\ge & {1 \over 6}(\frac{2k-6k^3+4k^4-k\sqrt{4k - 1}+4k^2\sqrt{4k - 1}}{(2k^2+1)^2})+{c_1 \over t}={\alpha_k \over 6}+{c_1 \over t}
\end{eqnarray}
where $t$ is a sufficiently large integer and $c_1=\frac{3k-18k^2+18k^3+24k^4+3k\sqrt{4k - 1}+6k^2\sqrt{4k - 1}-18k^3\sqrt{4k - 1}}{6(2k^2+1)^2} >0$.

Now suppose $\vec{y}=(y_1, y_2, ..., y_t)$ is an optimal vector of $\lambda(G^*(t))$. Let $n$ be large enough. By Remark \ref{remarkblow}, $3$-uniform graph $S_n=(\lfloor ny_1\rfloor, \ldots, \lfloor ny_{t}\rfloor)\otimes G^*(t)$ has density at least $\alpha_k+\frac{c_1}{2t}.$ Since $\pi({\cal F})\le \alpha_k$, some member $F$ of $\cal F$ is a subgraph of $S_n$ for $n$ sufficiently large. For such $F\in \cal F$, there exists a subgraph $M$ of $G^*(t)$ with $|V(M)|\le |V(F)|\leq s$ so that $F\subset (n, n, \ldots, n) \otimes M$. By Fact \ref{mono} and  Fact \ref{lblow}, we have
\begin{equation}\label{lambdasmall}
   \lambda(F) {{\rm Fact} \ \ref{mono} \atop \leq }\lambda ((n, n, \ldots, n) \otimes M){{\rm Fact} \ \ref{lblow} \atop =} \lambda (M).
\end{equation}

Theorem \ref{theo1} will follow  from the following  lemma to be proved in Section \ref{prooflemmaresult1}.

\begin{lemma}\label{lemmaresult1}
 Let $M$ be any subgraph of $G^*(t)$ with $|V(M)| \leq s$. Then
\begin{equation}
\lambda (M) \leq \frac{1}{6}\alpha_k
\end{equation}
holds.
\end{lemma}

Assuming that Lemma \ref{lemmaresult1} is true and applying Lemma \ref{lemmaresult1} to (\ref{lambdasmall}), we have $$\lambda(F) \le {1 \over 6}\alpha_k$$ which contradicts  our choice of $F$, i.e., contradicts that $\displaystyle \lambda(F) >{1 \over 6}\alpha_k$ for all $F \in \cal F$.  \q

\bigskip

To complete the proof of Theorem \ref{theo1}, what remains is to show Lemma \ref{lemmaresult1}.

\subsection{Proof of Lemma \ref{lemmaresult1}}\label{prooflemmaresult1}
By Fact \ref{mono}, we may assume that $M$ is an induced subgraph of $G^*(t)$. Let $$U_i=V(M)\cap V_i=\{v_1^i, v_2^i, \cdots, v_{s_i}^i\}.$$ So $s=s_1+\cdots+s_{2k+1}$.

Similar to Claim \ref{reduce0}, we have
\begin{claim}\label{reduce} If $N$ is the $3$-uniform graph formed from $M$ by removing the edges contained in $U_{2k+1}$ and inserting the edges $v^{2k+1}_{1}v^{2k+1}_2v^{2k+1}_j$, where $3\leq j \leq s_{2k+1}$, then $\lambda(M)\leq \lambda(N)$.
\end{claim}

By Claim \ref{reduce}, the proof of Lemma \ref{lemmaresult1} will be completed if we show that $\lambda(N)\leq {\alpha_k \over 6}$.
By Lemma \ref{symmetry}, there exists an optimal  vector $\vec{z}$ of $\lambda(N)$ such that
\begin{equation} \label{weights}
w(v_1^{2k+1})=w(v_2^{2k+1})\defeq\frac{a}{2}, \ \ w(v_3^{2k+1})=w(v_4^{2k+1})=\cdots =w(v^{2k+1}_{s_{2k+1}}) \defeq\frac{b}{s_{2k+1}-2},
\end{equation}
where $w(v)$ denotes the component of $\vec{z}$  corresponding to vertex $v$. Let $w_1$ be the sum of the components of $\vec{z}$ corresponding to all vertices in $\cup_{i=1}^{2k}U_i$. Then
\begin{eqnarray*}
&&\lambda(N) \le \bigg(\frac{w_1}{2k}\bigg)^3{2k\choose 3}+\bigg(\frac{w_1}{2k}\bigg)^2{2k\choose 2}(1-w_1)+w_1\bigg(\frac{a^2}{4}+ab+\frac{b^2}{2}\bigg)+\frac{a^2}{4}b,
\end{eqnarray*}
where $w_1+a+b=1.$

Note that if $b\le w_1$ or $a=0$ or $w_1\ge \frac{1}{2}$, then
\begin{eqnarray*}
\lambda(N)&\le&\bigg(\frac{w_1}{2k}\bigg)^3{2k\choose 3}+\bigg(\frac{w_1}{2k}\bigg)^2{2k\choose 2}(1-w_1)+w_1\bigg(\frac{a^2}{2}+ab+\frac{b^2}{2}\bigg)\\
&=&\bigg(\frac{w_1}{2k}\bigg)^3{2k\choose 3}+\bigg(\frac{w_1}{2k}\bigg)^2{2k\choose 2}(1-w_1)+w_1\frac{(1-w_1)^2}{2}\\
&\leq &\lim_{n\to\infty}\lambda(B(2k, n))={\alpha_k\over 6}.
\end{eqnarray*}
So we may always assume that $w_1<\frac{1}{2}$. Since $b=1-w_1-a$, then
\begin{eqnarray*}
\lambda(N)&\le&\bigg(\frac{w_1}{2k}\bigg)^3{2k\choose 3}+\bigg(\frac{w_1}{2k}\bigg)^2{2k\choose 2}(1-w_1)\\
&+&w_1\bigg(\frac{a^2}{4}+a(1-w_1-a)+\frac{(1-w_1-a)^2}{2}\bigg)+\frac{a^2}{4}(1-w_1-a)\triangleq f(a),
\end{eqnarray*}
where $w_1+a<1.$
\begin{eqnarray*}
f'(a)&=&w_1\bigg(\frac{a}{2}+(1-w_1-a)-a-(1-w_1-a)\bigg)+\frac{a}{2}(1-w_1-a)-\frac{a^2}{4}\\
&=&-\frac{3a^2}{4}+\frac{a}{2}-aw_1.
\end{eqnarray*}
Note that $w_1<\frac{1}{2}$, then $f(a)$ is increasing in $[0, \frac{2-4w_1}{3}]$ and decreasing in $[\frac{2-4w_1}{3}, 1]$.
So
\begin{eqnarray*}
f(a)&\le&f(\frac{2-4w_1}{3})\\
&=&\bigg(\frac{w_1}{2k}\bigg)^3{2k\choose 3}+\bigg(\frac{w_1}{2k}\bigg)^2{2k\choose 2}(1-w_1)+\frac{11w_1^3}{54}-\frac{5w_1^2}{9}+\frac{5w_1}{18}+\frac{1}{27}=g(w_1).
\end{eqnarray*}
Then $g'(w_1)=(\frac{1}{4k^2}-\frac{7}{18})w_1^2-(\frac{1}{2k}+\frac{1}{9})w_1+\frac{5}{18}.$ Solving $g'(w_1)=0$, we obtain that $w_1=\frac{\pm\sqrt{\frac{36}{81}+\frac{1}{9k}-\frac{1}{36k^2}}-\frac{1}{9}-\frac{1}{2k}}{\frac{7}{9}-\frac{1}{2k^2}}.$
Note that $\frac{-\sqrt{\frac{36}{81}+\frac{1}{9k}-\frac{1}{36k^2}}-\frac{1}{9}-\frac{1}{2k}}{\frac{7}{9}-\frac{1}{2k^2}}<0$ and $\frac{1}{4k^2}-\frac{7}{18}<0.$ We will show that $\frac{\sqrt{\frac{36}{81}+\frac{1}{9k}-\frac{1}{36k^2}}-\frac{1}{9}-\frac{1}{2k}}{\frac{7}{9}-\frac{1}{2k^2}}>\frac{1}{2}$. It's sufficient to show that $\sqrt{\frac{36}{81}+\frac{1}{9k}-\frac{1}{36k^2}}>\frac{1}{2}+\frac{1}{2k}-\frac{1}{4k^2}$. Note that $$\sqrt{\frac{36}{81}+\frac{1}{9k}-\frac{1}{36k^2}}>\frac{2}{3}>\frac{23}{36}\ge \frac{1}{2}+\frac{1}{2k}-\frac{1}{4k^2}$$
holds for $k\ge 3$. As for $k=2$, we have
$$\sqrt{\frac{36}{81}+\frac{1}{9k}-\frac{1}{36k^2}}=\frac{\sqrt{165}}{18}>\frac{11}{16}=\frac{1}{2}+\frac{1}{4}-\frac{1}{36}.$$
Since earlier discussion allows us to assume that $w_1<\frac{1}{2}$, therefore $g(w_1)$ is increasing in $[0, \frac{1}{2}]$. Note that
$$\frac{11w_1^3}{54}-\frac{5w_1^2}{9}+\frac{5w_1}{18}+\frac{1}{27}\bigg|_{w_1=\frac{1}{2}}=\frac{1}{16}=w_1\frac{(1-w_1)^2}{2}\bigg|_{w_1=\frac{1}{2}}.$$
Then
\begin{eqnarray*}
\lambda(N)&\le& g(\frac{1}{2})\\
&=&\bigg(\frac{w_1}{2k}\bigg)^3{2k\choose 3}+\bigg(\frac{w_1}{2k}\bigg)^2{2k\choose 2}(1-w_1)+\frac{11w_1^3}{54}-\frac{5w_1^2}{9}+\frac{5w_1}{18}+\frac{1}{27}\bigg|_{w_1=\frac{1}{2}}\\
&=&\bigg(\frac{w_1}{2k}\bigg)^3{2k\choose 3}+\bigg(\frac{w_1}{2k}\bigg)^2{2k\choose 2}(1-w_1)+w_1\frac{(1-w_1)^2}{2}\bigg|_{w_1=\frac{1}{2}}\\
&\leq& \lim_{n\to\infty}\lambda(B(2k, n))={\alpha_k\over 6}.
\end{eqnarray*}
\q

\end{document}